










\input amstex
\documentstyle{amsppt}
\magnification\magstep1
\NoRunningHeads
\vsize 8.0 truein
\define\k{\ref\key}
\define\pp{\pages}
\define\en{\endref}

\define\f{\frac}

\define\qq{\qquad}
\define\qqq{\quad\qquad}

\define\r1{\Bbb R^1}
\define\rr{\Bbb R}

\define\cb{\Cal P}

\define\p{\partial}

\define\nf{\infty}

\define\lp{ {L^p} }

\define\half{\f{1}{2}}

\define\ra{\rightarrow}
\define\al{\alpha}
\define\be{\beta}

\define\de{\delta}

\define\la{\lambda}

\hsize=30truecc
\baselineskip=16truept

\topmatter
\title A sharp estimate for the Hardy-Littlewood maximal function\endtitle
\author Loukas Grafakos$^*$, Stephen Montgomery-Smith$^{**}$, and
Olexei Motrunich$^{***}$\endauthor
\affil University of Missouri, Columbia and Princeton University\endaffil
\noindent
\address
Department of Mathematics, University of Missouri, Columbia, MO 65211,  \newline 
Department of Mathematics, University of Missouri, Columbia, MO 65211,  \newline 
Department of Physics, Princeton University, Princeton, NJ 08544
\endaddress
\email \newline
\noindent loukas\@math.missouri.edu,
stephen\@math.missouri.edu, 
 motrunich\@princeton.edu
\endemail
\thanks $\,\,\,\,\,\, ^*$Research partially supported by the University
of Missouri Research Board. \endthanks
\thanks $\,\,\, ^{**}$Research partially supported by the National
Science Foundation. \endthanks
\thanks $^{***}$Research partially supported by the University
of Missouri-Columbia Research Council. \endthanks
\abstract
The best constant in the usual $L^p$ norm inequality for the centered
Hardy-Littlewood maximal function on $\Bbb R^1$ is obtained for the
class of all ``peak-shaped'' functions. A function
on the line is called ``peak-shaped'' if it is positive and convex
except at one point.
The techniques we use include 
variational
methods.

AMS Classification (1991): 42B25
\endabstract
\endtopmatter

\document

{\bf 0. Introduction.}
\smallskip
Let
$$
(Mf)(x)= \sup_{\de >0} {1 \over 2 \de }\int_{x-\de}^{x+\de}
|f(t)| \, dt\tag0.1
$$
be the centered Hardy-Littlewood maximal operator on the
line.
This paper grew out from our attempt to find the operator norm
of $M$ on $\lp (\r1)$.

Since $M$ is a positive operator, we may restrict our attention
to positive functions.
 Let $\cb$ be the set of all 
positive  functions $f$ on  $\r1$, which are convex except at
one point (where we also allow $f$ to be discontinuous).
We call such functions ``peak-shaped.''

We were able to find the best constant in the inequality
$$
\|Mf\|_{L^p} \le C(p) \|f\|_{L^p}\qqq \text{for
$f$ in $\cb \cap \lp$.}\tag0.2
$$
for $1<p<\nf$.
It turns out that the best such $C(p)$ is the unique number $c_p$
which satisfies the equality
$$
M(|x|^{-1/p}) = c_p \ |x|^{-1/p}. \tag0.3
$$
Note that the function $|x|^{-1/p}$ is locally integrable, so
the left-hand side of (0.3) is well-defined. Strictly speaking, the
function $|x|^{-1/p}$ doesn't  belong to the space
$\cb \cap L^p$. It is, however, the pointwise limit of
a sequence of functions in  $\cb \cap L^p$ and the norm
ratio of this sequence converges to $c_p$.

Here is our main result.

\proclaim{Theorem}
The smallest possible constant $C(p)$ in the inequality
 $$
\|Mf\|_{L^p} \le C(p) \|f\|_{L^p}\qqq \text{for
f in $\cb \cap \lp$},
$$
is
$$
c_p=\sup_{\tau >1} {(\tau +1)^{p-1\over p} + (\tau -1)^{p-1\over p} \over
2\tau {p-1\over p}},\tag0.4
$$
that is, the constant in (0.3).
\endproclaim

One may ask the corresponding question when $p=1$. 
It was communicated to us by Jos\'e  Barrionuevo [Ba] that
the best constant $C_1$ in the weak type inequality
$$
|\{x: (Mf)(x)>\la \}|    \le C_1 { \|f\|_{L^1} \over \la }
$$
for $f$ in   $\cb \cap L^1$ is $C_1=1$. This result is sharp and is
analogous to our result when $p=1$. 
(In fact, this result is valid for the wider class of
positive functions that are increasing on
$(-\infty,c)$, and decreasing on $(c,\infty)$ for some number $c$.)

It is still a mystery what happens for general functions $f$.
It is conjectured in [DGS] that  $c_p$ is the operator norm of
the Hardy-Littlewood maximal function on $L^p(\Bbb R^1)$.
Our methods will not work for arbitrary functions and we  will
point out during the proof where they break down. 
For general functions $f \in L^1$, the conjecture used to be that
$C_1=3/2$.  However, it has recently been shown by Aldaz [Al]
that $C_1$ lies between $3/2$ and n    
$2$.  This result tends to suggest that the value $c_p$ given
by our Theorem is not the best constant for general $f\in L^p$.

The authors would like to thank the anonymous referee for many valuable
comments, and for pointing out errors in the original version of the
manuscript.

\bigskip

{\bf 1. Some preliminary Lemmas.}
\smallskip

Throughout this paper  we fix  a $p$ with $1<p<\nf$, and
a positive function $f$ in
$\cb \cap \lp$. 
Since $M$ commutes with translations, we may assume that
$f$ is convex except at $0$.  By a density argument, 
we may also make the following assumptions:
\roster
\item $f$ is smooth everywhere except at $0$;
\item $f$ has compact support;
\item $f(x) = |x|^{-1/2p}$ for $|x|$ sufficiently small.
\endroster
We remark that the third condition may seem unnatural, but 
it  simplifies some of the technicalities of the proof.
Equivalently, this condition may  be replaced by a  general 
assumption that $f$ is ``spiky'' enough near the origin but $f(0)<\infty$.

For every real $x\ne0$, define the function $\xi_x(t)=
\dfrac{1}{2t} \displaystyle\int_{x-t}^{x+t}f(u)\, du$ if $t>0$, and
$\xi_x(0) = f(x)$. It can be
seen that $\xi_x(t)$ is a $C^\nf$ function of $t> 0$
(except at $t = |x|$, where it is merely continuous)
and that it tends to zero
as $t\ra \nf$. 

Furthermore, we see that 
$$ \xi_x'(t) = {{1\over2}(f(x+t)+f(x-t)) - \xi_x(t) \over t} .$$
Convexity shows us that $\xi_x'(t) \ge 0$ for $t \in (0,|x|)$, and
the third condition on $f$ shows us that $\xi_x'(t) > 0$ for $t$ close
to $|x|$.  Thus we see that
$\xi_x(t)$ is non-decreasing for $t$ in some open neighborhood
of $(0, |x|]$.

The global maximum of $\xi_x$ over $[0,\nf)$ is
equal to $(Mf)(x)$. This maximum
is attained on some  set of  real numbers  $B_x=\{ t: \xi_x(t) =
\sup_{u\ge 0}\xi_x(u)\}$. Set $\de(x)= \max B_x$. Since $B_x$ is
a closed set, it contains $\de(x)$.  Note that $\de(x) > |x|$ for
$x\ne0$.
Thus
$$
(Mf)(x)= {1 \over 2 \de(x) }\int_{x-\de(x)}^{x+\de(x)}
f(t) \, dt.\tag1.1
$$
Since $\de(x)$ is a critical point of $\xi_x$, it follows that
$\xi_x'(\de(x))=0$. A simple calculation  and (1.1)
give formula (1.2) below.

Now fix $x_0 \ne 0$.
By the Implicit Function Theorem, the equation
$\xi_x'(\de)=0$ can be solved for $\de$ as a smooth
function of $x$ in the vicinity of any point $(x_0,\de(x_0))$,
as long as $\dfrac{ \p   \xi_x'(\de)}{ \p \de} \neq 0$ at $(x_0,\de(x_0))$.
This condition is equivalent to
$f'(x_0+\de(x_0))\neq f'(x_0-\de(x_0))$, which follows from the fact
that  $x_0+\de(x_0)$ and $x_0-\de(x_0)$ lie on opposite sides of
the origin and that
$f$ has different kind of monotonicity on each side.
Therefore $\de$ coincides with a smooth function in the
neighborhood of every point $x_0 \ne 0$, which implies that
$\de(x)$ is a smooth function of $x \ne 0$.
As a consequence $(Mf)(x)$ is also smooth for $x \ne 0$.

We notice that for sufficiently small $|x|$ that
$\de(x) = (1+\tau_{2p}) |x|$ for a fixed value $\tau_{2p}$, and
that $(Mf)(x) = c_{2p} |x|^{-1/2p}$.
Thus $\de(x)$ is a continuous function of $x$.

\proclaim{Lemma 1} For $x\ne 0$, we have
$$
 (Mf)(x)= {f(x+\de(x))+f(x-\de(x))\over 2},\tag1.2
 $$
 and
 $$
(Mf)'(x)= {f(x+\de(x))-f(x-\de(x))\over 2\de(x)}.\tag1.3
$$

\endproclaim

{\smc Proof.} $\,$ (1.2) is proved as indicated above.
To prove (1.3), differentiate the identity (1.1) and use (1.2).
This completes the proof of Lemma 1. $QED$.

\medskip

Formula (1.3) indicates that the points $x+\de(x)$ and
$x-\de(x)$ are  the $x$-coordinates of  some two points of
intersection of the graph of  $f$ with the
tangent line to the graph of $Mf$ at $(x,f(x))$.

\proclaim{Lemma 2} If $x>0$ then $\delta'(x) > 1$, and if
$x<0$ then $\delta'(x) < -1$.  Moreover $Mf$ is in $\cb$ with
its maximum at $0$.

\endproclaim

{\smc Proof.} $\,$
We begin by showing that $Mf$ has no inflection points away
from $0$. Differentiating (1.2) and (1.3) we obtain
that for $x\ne 0$, we have
$$
\align
 (Mf)'(x) &= f'(x+\de(x)){ (1+\de'(x))\over 2} +
  f'(x-\de(x)){ (1-\de'(x))\over 2 } \tag1.4 \\
 (Mf)'(x)\de'(x) +\de(x)(Mf)''(x) &=
f'(x+\de(x)){ (1+\de'(x))\over 2} -
  f'(x-\de(x)){ (1-\de'(x))\over 2 }. \tag1.5
\endalign
$$
If $q \ne 0$ were an inflection point,  then $(Mf)''(q)=0$,
and by (1.4) and (1.5) it would follow that
$$\align
f'(q+\de(q)) (1+\de'(q))&=(1+\de'(q))(Mf)'(q)\\
\noalign{\noindent or}
f'(q-\de(q)) (1-\de'(q))&=(1-\de'(q))(Mf)'(q).
\endalign$$
Then  $(Mf)'(q) $ would be   equal to either
$f'(q+\de(q))$ or $f'(q-\de(q))$. By
 Lemma 1,   $ (Mf)'(q)$ is the slope of the line
segment that joins $(q-\de(q),f(q-\de(q)))$ to
$(q+\de(q),f(q+\de(q)))$.  By the convexity conditions
on $f$, this line would then necessarily lie on 
the graph of $f$.  By (1.2), this would imply that
$(Mf)(q) \le f(q)$, a contradiction if condition~(3) is imposed
upon $f$.
Therefore $Mf$ has no inflection points away from $0$, hence  it is either concave or 
convex there.
Since $(Mf)(x)$ looks like ${1\over x}$ near $\pm \nf$, it follows that
$Mf$ is convex on $(-\infty,0)$ and on $(0,+\infty)$.

We now show that if $x < 0$, then $\de'(x) < -1$.
Let $x_1<x_2<0$ and let
 $L_i$ be the tangent line to the graph of $Mf$ at $x_i$,
$ i=1,2$. $L_i$ passes through the point
$(x_i+\de(x_i),f(x_i+\de(x_i)))$.
Since $Mf$ is convex on $(-\nf ,0)$, the
 line $L_1$ lies lower
than the
 line $L_2$ to the
right of $x_2$. Since $f$ is
decreasing on $(0,\nf)$, it follows that $L_2$ intersects
the graph of  $f$ on $(0,\nf)$ at a point with
$x$-coordinate bigger than the $x$-coordinate of the
intersection of $L_1$ with the graph of $f$.
This implies  that $x_1+\de(x_1)> x_2+\de(x_2) $ which
proves  that  $x+\de(x)$ is decreasing on $(-\nf ,a)$.
Therefore $\de'(x) < -1$ on $(-\nf ,0)$.
Likewise one can show that
 $\de'(x) > 1$ on $(0 ,+\nf)$. This completes the proof of Lemma 2. $QED$.

\bigskip

{\bf 2. The variational functional.}
\smallskip

It will be convenient to have $\de(x) \le 0$ for $x<0$.
To achieve this we define  $s(x)$ to be equal to  $\de(x)$  for $x>0$,
to be equal to $-\de(x)$ for $x<0$, and equal to $0$ if $x=0$.  
We  observe that
(1.1), (1.2), and (1.3) remain valid for $s(x)$.
We also observe that  $s'(x)>1$ for $x \ne 0$ and that
\roster

\item $\lim\limits_{x\ra +\nf}x+s(x)=  +\nf$,
\item $\lim\limits_{x\ra -\nf}x+s(x)=  -\nf$,
\item $\lim\limits_{x\ra +\nf}x-s(x)=  \inf (\text{support(f)})$,
\item $\lim\limits_{x\ra -\nf}x-s(x)=  \sup (\text{support(f)})$.

\endroster

For simplicity, we denote by $g=Mf$  the maximal function of $f$.
It turns out that a suitable convex combination of the integrals of
the functions $(g(x)-g'(x)s(x))^p(s'(x)-1)$ and
$(g(x)+g'(x)s(x))^p(s'(x)+1)$ will give rise to a
functional  related  to  $\|f\|_{L^p}^p$. Our goal will be to
minimize this functional by selecting a suitable $s(x)$. To
find such a minimizer,
we  solve the corresponding Euler-Langrange equations.

By Lemma 1 we have that
$$f(x+s(x))=g(x)+g'(x)s(x),\tag2.1$$
and
$$f(x-s(x))=g(x)-g'(x)s(x).\tag2.2$$
Raise both sides of (2.1) to the power $p$, multiply them
by $1+s'(x)$, and integrate from $-\nf$ to $\nf$ to obtain
$$
\align
&\int_{-\nf}^{+\nf} (g(x)+g'(x)  s(x))^p  (s'(x)+1)\,dx   \\
=& \int_{-\nf}^{+\nf} f(x+s(x))^p(s'(x)+1)\,dx  \\
=& \int_{-\nf}^{+\nf} f(x)^p\,dx     
= \|f\|_{L^p}^p   \tag2.3
 \endalign
 $$
Similarly, raise both sides of (2.2) to the power $p$,
multiply them by $s'(x)-1$, and integrate from $-\nf$ to $\nf$.
We obtain
$$
\align
&\int_{-\nf}^{+\nf} (g(x)-g'(x)  s(x))^p  (s'(x)-1)\,dx  \\
=& \int_{-\nf}^{+\nf} f(x-s(x))^p(s'(x)-1)\,dx    \\
=& \int_{-\nf}^{+\nf} f(-x)^p\,dx     
= \|f\|_{L^p}^p   \tag2.4
 \endalign
 $$

At this point, we remark that the calculation above will not work
for general functions $f$, because in that case the  function $s(x)$
will have many discontinuities, and above formulae will have
to include  terms needed to account for these discontinuities.
These discontinuities are generally rather unpredictable, and we have
not been able to find a way to deal with this problem.

Let $\half <\al <1$ be a real number to
 be selected later to depend on $p$ only.
Let $F$ be the following function of three variables:
$$
F(x,y,z)= \al (g(x)+g'(x)y)^p(z+1)
 + (1-\al) (g(x)-g'(x)y)^p(z-1).\tag2.5
$$
The domain of $F$ is the set of all $(x,y,z)$ which
satisfy
\roster
\item $-\nf <x<\nf$
\item $-\dfrac{g(x)}{ |g'(x)|}<y<\dfrac{g(x)}{ |g'(x)|}$
\item $ -\nf <z<\nf$
\endroster
Because of (2.1), (2.2) and the positivity of $f$, we have
that $(x,s(x),s'(x))$ lies in the domain of definition
of $F$.
Combining  (2.3) and (2.4) we obtain
$$
\|f\|_{L^p}^p = \int_{-\nf}^{+\nf} F(x,s(x),s'(x))\, dx.\tag2.6
$$

Now we define a functional  $I(\phi) $ that we would like
to minimize.  The domain of $I$ will be those functions
$\phi:\rr \to \rr$ that are smooth except possibly at $0$,
such that $F(x,\phi(x),\phi'(x))$ is integrable on compact subsets of 
$\rr \setminus \{0\}$, and such that the following improper integral 
converges:
$$
\leqalignno{
I(\phi) &= \int_{-\nf}^{+\nf} F(x,\phi(x),\phi'(x))\, dx \cr
&= \lim\Sb a\to-\nf\\ b\to0^-\\ c\to0^+\\ d\to+\nf\endSb
   \left(\int_a^b + \int_c^d \right) F(x,\phi(x),\phi'(x))\, dx .  \cr  }
$$
Denote by $\p_1 F$, $\p_2 F$, and $\p_3 F$ the
partial derivatives of $F(x,y,z)$ with respect to $x$, $y$, and $z$
respectively.
To minimize $I$,
we consider the associated Euler-Langrange
equations:
$$
{d  \over  dx}\big[(\p_3 F)(x,\phi(x),\phi'(x))\big] =
(\p_2 F)(x,\phi(x),\phi'(x)).\tag2.7
$$
We now have the following result.

\proclaim{Lemma 3} The function
$$ s_0(x)=-\be {g(x) \over g'(x)}, \qq \text{where}\qq \be= \be(\alpha) =
{ \al^{1\over p-1} - (1-\al)^{1\over p-1} \over
\al^{1\over p-1} +(1-\al)^{1\over p-1}  },\tag2.8
$$
is an exact solution of  equation (2.7) on $\Bbb R^1\setminus\{0\}$.
\endproclaim

\noindent {\smc Remark:}  Note that $0 < \be (\alpha) < 1$, since $\al > \half$.
\smallskip
 {\smc Proof.}
To prove Lemma 3, rewrite
 (2.7)   as
$$[ \al(g(x)+g'(x)\phi(x))^{p-1}-(1-\al) p (g(x)-g'(x)\phi(x))^{p-1}]
g''(x) \phi(x) =0.$$
Then substituting for $\phi = s_0$, we obtain the result. $QED.$

\smallskip

We would like to to be able to directly deduce that $I(s)\ge I(s_0)$. 
Unfortunately, general theorems from calculus of variations (see for example  [Br]) 
are not directly applicable here since $F(x,y,z)$ does not satisfy the 
usual convexity conditions needed. As it turns out, the desired inequality 
$I(s)\ge I(s_0)$ will be a consequence of the key inequality below which is 
true because of the very specific structure of the function $F(x,y,z)$. 

\proclaim{Lemma 4} For all $x\ne0$ we have
$$
\leqalignno{
 & F(x, s(x),s'(x)) -  F(x,s_0(x),s_0'(x))  \cr &\ge 
   (\p_2 F)(x,s_0(x),s_0'(x)) (s(x)-s_0(x)) \cr &
  \quad +
(\p_3 F)(x,s_0(x),s_0'(x)) (s'(x)-s_0'(x))  &(2.9)\cr }
$$
\endproclaim

{\smc Proof.}
We observe that if $h$ is a convex function on an interval
$J$, then for all $x,y$ in $J$ we have
$$
h(x)-h(y) \ge h'(y) (x-y)\tag2.10
$$
irrespectively of the order of $x$ and $y$. Next observe that
for all $x$ and all $z> 1$ the function $F(x,y,z)$ is
convex in $y$.  This is because when $z>1$,
$(\p_2^2 F)(x,y,z) >0$ for all $x,y$.

Let $x \ne 0$. Then by Lemma 2, $s'(x)=|\de'(x)|> 1$ and
$$
\leqalignno{
& F(x,s(x),s'(x))- F(x,s_0(x),s_0'(x)) \cr &= 
 [F(x,s(x),s'(x))-F(x,s_0(x),s'(x))] \cr &\quad +
[F(x,s_0(x),s'(x))-F(x,s_0(x),s_0'(x))] \cr &\ge 
    (\p_2 F)(x,s_0(x),s'(x)) (s(x)-s_0(x)) \cr &\quad +
(\p_3 F)(x,s_0(x),s_0'(x)) (s'(x)-s_0'(x)),   &(2.11)   \cr }
$$
by convexity of $F$ in the second variable, (2.10), and linearity
of $F$ in the third variable.

Calculation gives
$$
(\p_3\p_2F)(x,y,z)=pg'(x)\big[\al (g(x)+g'(x)y)^{p-1} -
(1-\al) (g(x)-g'(x)y)^{p-1}\big].\tag2.12
$$
Setting $y=s_0(x)$ in (2.12) we obtain that
$$
(\p_3\p_2F)(x,s_0(x),z)=pg'(x)\big[\al (1-\be)^{p-1} -
(1-\al) (1+\be)^{p-1}\big]=0,
$$
since by the definition of $\be$, it follows that
$$
\al = { (1+\be)^{p-1} \over  (1+\be)^{p-1} + (1-\be)^{p-1}}.
$$
We have now proved that the function $(\p_2F)(x,s_0(x),z)$
is constant in $z$.  The proof of (2.9) is now complete
if we replace $s'(x)$  by $s_0'(x)$ in the first summand of (2.11).
$QED.$

\bigskip

{\bf 3. The core of the proof.}
\smallskip

Next we have the following. 

\proclaim{Lemma 5}
Both $s$ and $s_0$ lie in the domain of the functional
$I$.  We have the equality
$$ I(s_0) = r(\alpha) \|g\|_{L^p}^p ,\tag 3.1$$
where $r(\alpha) = \gamma_1(\alpha) + p \beta(\alpha) \gamma_2 (\alpha)$
with $\gamma_1(\alpha) = \alpha(1-\beta)^p - (1-\alpha)(1+\beta)^p$ and
$\gamma_2(\alpha) = \alpha(1-\beta(\alpha))^p + (1-\alpha)(1+\beta(\alpha))^p$.  We also
have the inequality
$$ I(s) = \|f\|_{L^p}^p \ge I(s_0). \tag 3.2$$
\endproclaim

{\smc Proof.}
First, it is clear that $s$ is in the domain of $I$, by the calculations
at the beginning of the previous section.  Let us work with $s_0$.  We see
that if $0<a<b<\infty$, then integrating by
parts we get
$$ \eqalignno{
   \int_a^b F(x,s_0(x),s_0'(x)) \, dx
   &=
   \int_a^b \gamma_1 g(x)^p - \gamma_2 g(x)^p {d\over dx}
   \left( \beta {g(x) \over g'(x)} \right) \, dx \cr
   &=
   r(\alpha) \int_a^b g(x)^p \, dx
   + \gamma_2\beta
   \left( {g(a)^{p+1} \over g'(a)} - {g(b)^{p+1} \over g'(b) }\right) . \cr } $$
As $a \to 0$, we have the explicit formula for $g(a) = c_{2p} |a|^{-1/2p}$
which tells us that ${g(a)^{p+1} \over g'(a) } \to 0$.  When $b$ is very
large, we note that $f(b+\delta(b)) = 0$ (because $f$ has compact support),
and hence by (1.2) and (1.3)
we have that $\left|{g(b) \over g'(b)}\right| = \delta(b)$.  
Furthermore,
$|b+\inf (\text{support(f)})| \le \delta(b) \le |b+\sup (\text{support(f)})|$,
and $g(b) = O\left({1\over b}\right)$. Thus ${g(b)^{p+1} \over g'(b)}
\to 0$ as $b \to +\nf$.  We obtain a similar result if $-\infty <a<b<0$.
Since $0<\beta<1$, we obtain that 
$s_0$ is in the
domain of $I$, and that $ I(s_0) = r(\alpha) \|g\|_{L^p}^p $.

Now let us consider $I(s)$, which we already  know is equal to 
$\|f\|_{L^p}^p$.
It is here that estimate (2.9) plays its crucial role.
If $0<a<b<\infty$, then
$$
\leqalignno{
 \int_{a}^{b} &F(x,s(x),s'(x))\, dx \cr
&= \int_{a}^{b} [F(x,s(x),s'(x))-F(x,s_0(x),s_0'(x))]\, dx 
   + \int_{a}^{b} F(x,s_0(x),s_0'(x))\, dx   \cr
 & \ge  \int_{a}^{b} [(\p_2F)(x,s_0(x),s_0'(x))(s(x)-s_0(x)) +
 (\p_3F)(x,s_0(x),s_0'(x))(s'(x)-s_0'(x))]\, dx   &\cr
 &\qqq\qqq\qqq\qqq +\int_{a}^{b} F(x,s_0(x),s_0'(x))\, dx,  &(3.3)\cr
}
$$
where we used Lemma 4 in the inequality above.
Next, we integrate by parts, and
(3.3) is now equal to
$$
\leqalignno{
&\qq \int_{a}^{b}     \bigg[ (\p_2 F)(x,s_0(x),s_0'(x)) -
{d \over dx}\big( (\p_3 F)(x, s_0(x), s_0'(x)) \big) \bigg]
(s(x)-s_0(x))\, dx           \cr
&\qqq + (\p_3 F)(b,s_0(b),s_0'(b)) (s(b)-s_0(b)) 
      - (\p_3 F)(a,s_0(a),s_0'(a)) (s(a)-s_0(a))  \cr
&\qqq + \int_{a}^{b} F(x,s_0(x),s_0'(x))\, dx.    \cr
}
$$
First note that the first integral evaluates to $0$ by
Lemma 3.  Now,
as $a \to 0$, we may explicitly calculate to see that
$(\p_3 F)(a,s_0(a),s_0'(a)) (s(a)-s_0(a)) \to 0$.  Also,
if $b$ is very large, using the fact that
$f(b+\delta(b)) = 0$, (1.2), (1.3), and (2.8), we see that $s_0(b) = \beta s(b)$,
and hence $(\p_3 F)(b,s_0(b),s_0'(b)) (s(b)-s_0(b))$ is 
a constant multiple of $g(b)^p s(b)$.  Arguing as above,
we see that $ (\p_3 F)(b,s_0(b),s_0'(b)) (s(b)-s_0(b)) $
 tends to $0$ as $b \to \infty$. We obtain a similar result 
when $-\infty <a<b<0$, hence we conclude that 
$I(s)\ge I(s_0)$. 
$QED.$

\bigskip

To finish the proof, we only need the following result.

\proclaim {Lemma 6}
There exists ${1\over2} < \alpha < 1$ such that
$$ r(\alpha) = c_p^{-p} ,$$
where $c_p$ is the constant in (0.4). 
\endproclaim

 \noindent {\smc Remark:}
In fact it is true that $c^{-p}_p$ is the absolute minimum of $r(\alpha)$
for ${1\over2} < \alpha < 1$. 
\smallskip
{\smc Proof.}
We have that
$$
r(\alpha)=\alpha (1-\beta)^p - (1-\alpha)(1+\beta)^p +p \beta \gamma_2 =
{  2^p (p-1)\alpha (1-\alpha)
(\alpha^{1 \over p-1} -(1-\alpha)^{1 \over p-1})
  \over (\alpha^{1 \over p-1} +(1-\alpha)^{1 \over p-1})^p }.
$$
For $t\in (1,+\nf)$, define
$$
h(t)={ (t+1)^{p-1 \over p} + (t-1)^{p-1 \over p}
\over 2\, {p-1\over p}\, t}.
$$
We see that $h'(t)$ is monotonically decreasing on $(1,+\infty)$ and
attains its only zero  at the unique $\tau$ satisfying
$$\qquad \bigg({p+\tau \over p-\tau}\bigg)^p =
{\tau+1 \over \tau-1}. \tag3.4
$$
Let
$$\alpha_0 = { (p+\tau)^{p-1}  \over  (p+\tau)^{p-1} +
(p-\tau)^{p-1}}\tag3.5
$$
where $\tau$ satisfies (3.4). It is clear that $\half < \al_0 <1$.
Then by elementary, although perhaps not easy manipulations, we have
$$
\leqalignno{
r(\alpha_0) = &
{ 2\, (p-1)\,  \tau \, (p+\tau)^{p-1} (p-\tau)^{p-1} \over
p^p\, [(p+\tau)^{p-1}+(p-\tau)^{p-1}]}  \cr
= &  {2^p \, (p-1)\,  \tau \over p \, ({p-\tau \over p+\tau} +1)^{p-1}
(({p+\tau \over p-\tau})^{p-1}+1)} \cr
= &{2^p\,  (p-1)\,  \tau  \over
p\, [(\tau-1)^{-{1\over p}}+(\tau+1)^{-{1\over p}}]^{p-1}
[(\tau+1)^{p-1 \over p}+(\tau-1)^{p-1 \over p}]}\cr
= & \bigg( { 2\tau \,{p-1\over p}   \over  (\tau-1)^{{p-1\over p}}+
(\tau+1)^{{p-1\over p}} } \bigg)^p = h(\tau)^{-p}
= c_p^{-p}, &(3.6) \cr
}
$$
where $c_p$ is the constant in  (0.4).
$QED.$

\enddocument